\input amssym.def
\magnification 1200
\overfullrule0pt
\hyphenation{pre-print}
\nopagenumbers
\headline={\ifnum\pageno=1 \hfill \else\hss{\tenrm--\folio--}\hss \fi}
\newcount\sectionnumber
\newcount\equationnumber
\newcount\refnumber
\def\ifundefined#1{\expandafter\ifx\csname#1\endcsname\relax}
\def\assignnumber#1#2{%
	\ifundefined{#1}\relax\else\message{#1 already defined}\fi
	\expandafter\xdef\csname#1\endcsname
 {\the\sectionnumber.\the#2}}%
\def\secassignnumber#1#2{%
  \ifundefined{#1}\relax\else\message{#1 already defined}\fi
  \expandafter\xdef\csname#1\endcsname{\the#2}}%
%
%
\def\secname#1{\relax
  \global\advance\sectionnumber by 1
  \secassignnumber{S#1}\sectionnumber
  \csname S#1\endcsname}
\def\Sec#1 #2 {\vskip0pt plus.1\vsize\penalty-250\vskip0pt plus-.1\vsize
  \bigbreak\bigskip
  \equationnumber0
  \noindent{\bf \secname{#1}. #2}\par
  \nobreak\smallskip\noindent}
\def\sectag#1{\ifundefined{S#1}\message{S#1 undefined}{\sl #1}%
  \else\csname S#1\endcsname\fi}
%
%
\def\eq#1{\relax
  \global\advance\equationnumber by 1
  \assignnumber{EN#1}\equationnumber
  {\rm (\csname EN#1\endcsname)}}
\def\eqtag#1{\ifundefined{EN#1}\message{EN#1 undefined}{\sl (#1)}%
  \else(\csname EN#1\endcsname)\fi}
%
%
\def\refitem#1 #2\par{\ifundefined{REF#1}
\global\advance\refnumber by1%
\expandafter\xdef\csname REF#1\endcsname{\the\refnumber}%
\else\item{\ref{#1}}#2\sLP\fi}

\def\ref#1{\ifundefined{REF#1}\message{REF#1 is undefined}\else
  [\csname REF#1\endcsname]\fi}
\def\reff#1#2{\ifundefined{REF#1}\message{REF#1 is undefined}\else
  [\csname REF#1\endcsname, #2]\fi}
\def\Ref{\vskip0pt plus.1\vsize\penalty-250\vskip0pt plus-.1\vsize
  \bigbreak\bigskip\leftline{\bf References}\nobreak\smallskip
  \frenchspacing}
%
%

\def\Proof{\sLP{\bf Proof}\quad}
\def\halmos{\hbox{\vrule height0.31cm width0.01cm\vbox{\hrule height
 0.01cm width0.3cm \vskip0.29cm \hrule height 0.01cm width0.3cm}\vrule
 height0.31cm width 0.01cm}}
\def\hhalmos{{\unskip\nobreak\hfil\penalty50
	\quad\vadjust{}\nobreak\hfil\halmos
	\parfillskip=0pt\finalhyphendemerits=0\par}}

\let\mPP=\medbreak
\let\bPP=\bigbreak
\def\LP{\par\noindent}
\def\sLP{\smallbreak\noindent}
\def\mLP{\medbreak\noindent}
\def\bLP{\bigbreak\noindent}
\let\al\alpha
\let\de\delta
\let\la\lambda
\let\si\sigma
\let\De\Delta
\def\CC{{\Bbb C}}
\def\JJ{{\Bbb J}}
\def\QQ{{\Bbb Q}}
\def\RR{{\Bbb R}}
\def\FSD{{\cal D}}
\def\FSF{{\cal F}}
\def\FSR{{\cal R}}
\def\gog{{\goth g}}
\def\goh{{\goth h}}
\def\gok{{\goth k}}
\def\gon{{\goth n}}
\let\lan=\langle
\let\ran=\rangle
\let\ten=\otimes
\let\dirsum=\oplus
\def\id{{\rm id}}

\def\wt{{\rm wt}}
\def\Ind{{\rm Ind}}
\def\Hom{{\rm Hom}}
\def\End{{\rm End}}
\def\Tr{{\rm Tr}}
\let\lrarrow\leftrightarrow

\font\titlefont=cmr10 at 17.28truept
\font\authorfont=cmti10 at 14truept
\font\addressfont=cmr10 at 10truept
\font\ttaddressfont=cmtt10 at 10truept

\refitem B-B-B

\refitem E-S

\refitem E-V1

\refitem E-V2

\refitem E-V3

\refitem E-V4

\refitem Fe1

\refitem Fe2

\refitem G-N

\LP
{\titlefont Some details of proofs of theorems related to the
\mLP
quantum dynamical Yang-Baxter equation}
\bLP
{\authorfont Tom H. Koornwinder}
\bLP
Version of July 13, 2000

\bLP
{\bf Abstract}
\LP
This paper gives some further details of proofs of some theorems
related to the quantum dynamical Yang-Baxter equation.
This mainly expands proofs given in ``Lectures on the
dynamical Yang-Baxter equation'' by P. Etingof and O. Schiffmann,
math.QA/9908064. This concerns the intertwining operator, the fusion
matrix, the exchange matrix and the difference operators.
The last part expands proofs
given in ``Traces of intertwiners for quantum groups and difference
equations, I'' by P. Etingof and A. Varchenko, math.QA/9907181.
This concerns the dual Macdonald-Ruijsenaars equations.
This paper does not claim originality, priority or completeness.
It is meant as a service to whoever may take profit of it.

\Sec{intro} {Introduction}
The quantum dynamical Yang-Baxter equation (QDYBE)
was first considered in 1984
by Gervais and Neveu \ref{G-N}, with motivation from physics
(for monodromy matrices in Liouville theory).
A general form of QDYBE with spectral parameter was presented by
Felder \ref{Fe1}, \ref{Fe2} at two major congresses in 1994.
The corresponding classical dynamical Yang-Baxter equation (CDYBE)
was presented there as well.
Next Etingof and Varchenko started a program to give geometric interpretations
of solutions of CDYBE (see \ref{E-V1}) and of QDYBE (see \ref{E-V2})
in the case without spectral parameter. In the context of this
program they pointed out a method to obtain solutions of QDYBE by
the so-called exchange construction (see \ref{E-V3}). This uses,
for any simple Lie algebra $\gog$,
representation theory of $U(\gog)$ or of its quantized version $U_q(\gog)$
in order to define successively the intertwining operator, the fusion matrix
and the exchange matrix. The matrix elements of the intertwining operator
and of the exchange matrix generalize respectively the Clebsch-Gordan
coefficients and the Racah coefficients to the case where the
first tensor factor is a Verma module rather than a finite dimensional
irreducible module. The exchange matrix is shown to satisfy QDYBE.
Etingof and Varchenko also started in \ref{E-V4}
a related program to connect the above
objects with weighted trace functions and with solutions of the
(q-)Knizhnik-Zamolodchikov-Bernard equation (KZB or qKZB).

A nice introduction to the topics indicated above was recently given by
Etingof and O. Schiffmann \ref{E-S}. While I was reading this paper
in connection with a seminar in Amsterdam during the fall of 1999,
I added some details of proofs for my own convenience, and I put these
notes in TeX in order that the other participants in the seminar
could take profit of it. I put these informal notes on my homepage.
Since the version v2 of \ref{E-S} is now referring to these notes,
I decided to post the paper on QA.

I want to emphasize that these notes are purely meant as a service
to whoever may take profit of it. I do not claim any originality or
priority with these proofs. Neither I tried to cover the full contents
of \ref{E-S}. Most of my paper only treats the $q=1$ case.
Only the second part of the section on the exchange matrix also
covers the quantum case.
In general, the extension to the quantum case will ususally be straightforward.

As for the contents, Sections 2, 3 and 4 respectively deal with the
intertwining operator, the fusion matrix and the exchange matrix.
In \ref{E-S} these topics are all covered in Section 2. My Sections 5 on
difference operators and 6 on weighted trace functions address some topics
in Section 9 of \ref{E-S} (Transfer matrices and generalized
Macdonald-Ruijsenaars equations). The details of proofs in Section 6
concern $q=1$ analogues of proofs given in Section 3 of \ref{E-V4} in
connection with the dual Macdonald-Ruijsenaars equations.

I want to call attention to one conceptual aspect.
This concerns formulas (4.8), (4.9). The first formula expresses an
exchange matrix $R_{U,V\ten W}(\la)$ after shifted conjugation by
the fusion matrix $J_{VW}(\la)$ as a product of $R_{UV}(\la)$
(with appropriately shifted $\la$) and $R_{UW}(\la)$. The second formula
is analogous. These formulas are not explicitly given in \ref{E-S},
but they do occur in \ref{E-V4} without getting particular emphasis.
They can be used in order to prove that $R(\la)$ satisfies QDYBE.
This is analogous to the role of the quasi-triangularity property
of the (non-dynamical) universal $R$-matrix for proving the QYBE in that case.
In fact, it is possible to see (4.8) and (4.9) in the context
of a certain quasitriangular quasi-Hopf algebra, see Babelon, Bernard \& Biley
\reff{B-B-B}{Section 3} for the quantum $sl(2)$ case.

\mLP
{\sl Acknowledgements}\quad
I was inspired by Pavel Etingof's lectures on the dynamical
Yang-Baxter equations at the London Mathematical Society Symposium
on Quantum groups in Durham, UK, July 1999.
\LP
I thank Eric Opdam for suggesting a shorter proof than I originally had
for the rational dependence on $\la$ of the intertwining operator.
\mLP
{\sl Notation}\quad
Throughout this paper I will denote by [E-S] the paper \ref{E-S}
by Etingof \& Schiffmann, and by [E-V] the paper \ref{E-V4} by Etingof
\& Varchenko.

\Sec{intertwining} {The intertwining operator}
First I make two preliminary remarks in preparation of the proof
of [E-S], Proposition 2.2.

\mPP
Let $\gog$ be a Lie algebra with Lie subalgebra $\gok$,
and let $V$ be a $\gok$-module. Then:
\sLP
$\Ind_\gok^\gog V:=U(\gog)\ten_\gok V$\quad with\quad
$a\cdot(u\ten_\gok v):=(au)\ten_\gok v$\quad
($a\in\gog$, $u\in U(\gog)$, $v\in V$).
\sLP
Let $W$ be a $\gog$-module. Then {\sl Frobenius reciprocity} states
that there is an isomorphism of linear spaces
$$
f\lrarrow F\colon \Hom_\gok(V,W)\lrarrow
\Hom_\gog(U(\gog)\ten_\gok V,W)
$$
given by\quad $F(u\ten_\gok v):=u\cdot f(v)$,\quad
$f(v):=F(1\ten_\gok v)$\quad($u\in U(\gog)$, $v\in V$).

\mPP
For the other remark let $\gog$ be a Lie algebra and let $Z,W,V$ be
$\gog$-modules. Then there is an isomorphism of linear spaces
$$
f\lrarrow F\colon\Hom_\gog(Z,W\ten V)\lrarrow
\Hom_\gog(Z\ten W^*,V)
$$
given by\quad$F(z\ten w^*)=\lan f(z),w^*\ran$\quad
($z\in Z$, $w^*\in W^*$).

\bLP
{\bf Proof of [E-S], Proposition 2.2.}\quad
We have a composition of five isomorphisms
$$\eqalignno{
&\Phi\lrarrow\Phi_1\lrarrow\Phi_2\lrarrow\Phi_3\lrarrow\Phi_4\lrarrow\Phi_5
=\lan\Phi\ran\colon\quad
\Hom_\gog(U(\gog)\ten_{\goh\ten \gon_+}\CC_\la,M_\mu\ten V)\lrarrow
\cr
&\quad\lrarrow\Hom_{\goh\ten \gon_+}(\CC_\la,M_\mu\ten V)
\lrarrow\Hom_{\goh\ten \gon_+}(\CC_\la\ten M_\mu^*,V)\lrarrow
\cr
&\qquad\lrarrow
\Hom_{\goh\ten \gon_+}(U(\gon_+)\ten_\goh\CC_{-\mu},V\ten\CC_\la^*)
\lrarrow\Hom_\goh(\CC_{-\mu},V\ten\CC_\la^*)
\lrarrow\Hom_\goh(\CC_\la\ten\CC_{-\mu},V),
\cr}
$$
where\quad
$\Phi_1(x_\la):=\Phi(x_\la)$,\quad
$\Phi_2(x_\la\ten u^*):=\lan\Phi(x_\la),u^*\ran$\quad($u^*\in M_\mu^*$),
\sLP
$\Phi_3(u^*):=\lan\Phi(x_\la),u^*\ran\ten x_\la^*$\quad
($u^*\in M_\mu^*\simeq\CC_{-\mu}\ten_\goh U(\gon_+)$),\quad
$\Phi_4(x_{-\mu}):=\lan\Phi(x_\la),x_\mu^*\ran\ten x_\la$,
\sLP
$\Phi_5(x_\la\ten x_{-\mu}):=\lan\Phi(x_\la),x_\mu^*\ran
=\lan\Phi\ran$.\hhalmos

\bLP
{\bf Proof that the coefficients of $\Phi_\la^v$ are rational in $\la$}\quad
(statement in paragraph after the proof of [E-S], Proposition 2.2;
the proof below is essentially due to Eric Opdam)
\sLP
Let $\al_1,\ldots,\al_N$ be the positive roots (the elements of $\De^+$).
Let $V$ be a finite-dimensional $\gog$-module, and let $v\in V\backslash\{0\}$
be $\goh$-homogeneous.
Consider the Verma module $M_{\la-\wt(v)}$ for generic values of
$\la\in\goh^*$, where it is irreducible.
By Proposition 2.2 there is a unique $\gog$-intertwining linear map
$\Phi_\la^v\colon M_\la\to M_{\la-\wt(v)}\ten V$ such that
$$
\Phi_\la^v(x_\la)=\sum_{k_1,\ldots,k_N\ge 0}
f_{\al_1}^{k_1}\ldots f_{\al_N}^{k_N}\cdot x_\mu\ten v_{k_1,\ldots,k_N}
\quad\hbox{with $v_{0,\ldots,0}=v$.}
\eqno\eq{1}
$$
Here $\mu:=\la-\wt(v)$. Clearly
$\wt(v_{k_1,\ldots,k_N})=\la-\mu+k_1\al_1+\cdots+k_N\al_N$.
It is sufficient to show that the $v_{k_1,\ldots,k_N}$ are rational in $\la$.

The unique existence of $\Phi_\la^v$ satisfying the above conditions
is equivalent to the unique existence of $w\in M_\mu\ten V$
such that $\wt(w)=\la$, $e_{\al_i}\cdot w=0$ for $i=1,\ldots,N$ and
such that
$w$ has the form of the right-hand side of \eqtag{1} with
$v_{0,\ldots,0}=v$.
We will show that the unique existence of $w$ with these properties
implies that the $v_{k_1,\ldots,k_N}$ are rational in $\la$.

Note that
$$
e_{\al_i}\,f_{\al_1}^{k_1}\ldots f_{\al_N}^{k_N}\cdot x_\mu
=\sum_{\matrix{\scriptstyle l_1,\ldots,l_N\ge0
\cr
\scriptstyle k_1\al_1+\cdots+k_N\al_N=
\cr
\scriptstyle \al_i+
l_1\al_1+\cdots+l_N\al_N}}
p_{i;l_1,\ldots,l_N}^{k_1,\ldots,k_N}(\la)\,
f_{\al_1}^{l_1}\ldots f_{\al_N}^{l_N}\cdot x_\mu
$$
with
$p_{i;l_1,\ldots,l_N}^{k_1,\ldots,k_N}(\la)$ polynomial in $\la$.
So, for $i=1,\ldots,N$ we have
$$
0=e_{\al_i}\cdot w=
\sum_{l_1,\ldots,l_N} f_{\al_1}^{l_1}\ldots f_{\al_N}^{l_N}
\cdot x_\mu\ten
\Bigl(e_{\al_i}\cdot v_{l_1,\ldots,v_N}+
\sum_{\matrix{\scriptstyle k_1,\ldots,k_N\ge0
\cr
\scriptstyle k_1\al_1+\cdots+k_N\al_N=
\cr
\scriptstyle \al_i+
l_1\al_1+\cdots+l_N\al_N}}
p_{i;l_1,\ldots,l_N}^{k_1,\ldots,k_N}(\la)\,
v_{k_1,\ldots,k_N}\Bigr).
$$
So the inhomogeneous system of linear equations in
the coordinates of the vectors $v_{l_1,\ldots,l_N}$ ($l_1,\ldots, l_N$
nonnegative integers, not all 0) given by
$$
e_{\al_i}\cdot v_{l_1,\ldots,v_N}+
\sum_{\matrix{\scriptstyle k_1,\ldots,k_N\ge0
\cr
\scriptstyle k_1\al_1+\cdots+k_N\al_N=
\cr
\scriptstyle \al_i+
l_1\al_1+\cdots+l_N\al_N}}
p_{i;l_1,\ldots,l_N}^{k_1,\ldots,k_N}(\la)\,
v_{k_1,\ldots,k_N}=0\quad(i=1,\ldots,N)
$$
has for generic $\la$ a unique
solution. Since the coefficients are polynomials in $\la$ it follows
that the solution must be rational in $\la$.\hhalmos

\Sec{fusion} {The fusion matrix}
{\bf Proof of [E-S], Proposition 2.3, part 2}
\sLP
$$
\Phi_\la^v(x_\la)\in x_{\la-\wt(v)}\ten v+
M_{\la-\wt(v)}[<\la-\wt(v)]\ten V[>\wt(v)].
$$
Hence
$$
\Phi_\la^v(M_\la[<\la])\i x_{\la-\wt(v)}\ten V[<\wt(v)]+
M_{\la-\wt(v)}[<\la-\wt(v)]\ten V.
$$
It follows that
$$
\eqalignno{
&(\Phi_{\la-\wt(v)}^w\ten 1)(\Phi_\la^v(x_\la))
\cr
&\qquad\in \Phi_{\la-\wt(v)}^w(x_{\la-\wt(v)})\ten v+
\Phi_{\la-\wt(v)}^w(M_{\la-\wt(v)}[<\la-\wt(v)])\ten V[>\wt(v)]
\cr
&\qquad\i x_{\la-\wt(v)-\wt(w)}\ten w\ten v
+M_{\la-\wt(v)-\wt(w)}[<\la-\wt(v)-\wt(w)]\ten W\ten V
\cr
&\qquad\qquad+x_{\la-\wt(v)-\wt(w)}\ten W[<\wt(w)]\ten V[>\wt(v)]
\cr
&\qquad\qquad+M_{\la-\wt(v)-\wt(w)}[<\la-\wt(v)-\wt(w)]\ten W\ten V.
}
$$
Hence
$$
J_{WV}(\la)(w\ten v)\in
w\ten v+W[<\wt(w)]\ten V[>\wt(v)].
\eqno\halmos
$$

\bLP
{\bf Proof of [E-S], Proposition 2.3, part 3}\quad
On the one hand we have
$$
\eqalignno{
&(\Phi_{\la-\wt(v)-\wt(w)}^u\ten 1\ten 1)\circ
(\Phi_{\la-\wt(v)}^w\ten 1)\circ \Phi_\la^v(x_\la)&\eq{9}
\cr
&\quad=(\Phi_{\la-\wt(v)-\wt(w)}^u\ten 1\ten 1)\circ
\Phi_\la^{J_{WV}(\la)(w\ten v)}(x_\la)
\cr
&\quad=\Phi_\la^{J_{U,W\ten V}(\la)\circ(1\ten J_{WV}(\la))(u\ten w\ten v)}
(x_\la).&\eq{10}
}
$$
On the other hand, expression \eqtag{9} also equals
$$
\eqalignno{
&(\Phi_{\la-\wt(v)}^{J_{UW}(\la-\wt(v))(u\ten w)}\ten 1)
\circ \Phi_\la^v(x_\la)
\cr
&=\Phi_\la^{(J_{U\ten W,V}\ten 1)(\la)\circ J_{UW}(\la-\wt(v))(u\ten w\ten v)}
(x_\la).&\eq{11}
}
$$
Hence, by equality of expressions \eqtag{10} and \eqtag{11}, we have
$$
J_{U,W\ten V}(\la)\circ(1\ten J_{WV}(\la))(u\ten w\ten v)=
(J_{U\ten W,V}\ten 1)(\la)\circ J_{UW}(\la-\wt(v))(u\ten w\ten v).
$$
Hence we arrive at the {\sl dynamical 2-cocycle condition},
which was to be proved:
$$
J_{U,W\ten V}(\la)\circ(1\ten J_{WV}(\la))=
(J_{U\ten W,V}\ten 1)(\la)\circ J_{UW}(\la-h^{(3)}).
\eqno\eq{34}\;\;\halmos
$$

\Sec{exchange} {The exchange matrix}
Proposition 2.4 in [E-S] states that the exchange matrix
$R_{VW}(\la):=J_{VW}(\la)^{-1}\,J_{WV}^{21}(\la)$
satisfies the QDYBE
$$
R_{VW}(\la-h^{(3)})\,R_{VU}(\la)\,R_{WU}(\la-h^{(1)})=
R_{WU}(\la)\,R_{VU}(\la-h^{(2)})\,R_{VW}(\la)
\eqno\eq{39}
$$
as an identity of operators on $V\ten W\ten U$.

In preparation of the proof
recall that
$\Phi_\la^{w,v}:=(\Phi_{\la-\wt(v)}^w\ten 1)\circ\Phi_\la^v\,$.
Then [E-S] state:
\bLP
{\bf Lemma}\quad
$R_{VW}(\la)(v\ten w)=\sum_i v_i\ten w_i$\quad where\quad
$\Phi_\la^{w,v}=(1\ten P)\sum_i\Phi_\la^{v_i,w_i}$.
\mLP
{\bf \Proof}\quad
Assume $R_{VW}(\la)(v\ten w)=\sum_i v_i\ten w_i$. Then
$$
\eqalignno{
\Phi_\la^{w,v}
=\Phi_\la^{J_{WV}(\la)(w\ten v)}
=\Phi_\la^{PJ_{VW}(\la)R_{VW}(\la)(v\ten w)}
=(1\ten P)\Phi_\la^{J_{VW}(\la)R_{VW}(\la)(v\ten w)}\qquad&
\cr
=(1\ten P)\sum_i\Phi_\la^{J_{VW}(\la)(v_i\ten w_i)}
=(1\ten P)\sum_i\Phi_\la^{v_i,w_i}.&&\halmos
\cr}
$$
\bLP
{\bf First proof of QDYBE \eqtag{39}}\quad
Put
$$
\eqalignno{
\Phi_\la^{u,w,v}
&:=(\Phi_{\la-\wt(v)-\wt(w)}^u\ten 1\ten 1)\circ(\Phi_{\la-\wt(v)}^w\ten 1)
\circ\Phi_\la^v
\cr
&=(\Phi_{\la-\wt(v)}^{u,w}\ten 1)\circ\Phi_\la^v
\cr
&=(\Phi_{\la-\wt(v)-\wt(w)}^u\ten 1\ten 1)\circ\Phi_\la^{w,v}.
\cr}
$$
Now we have on the one hand
$$
\eqalignno{
\Phi_\la^{u,w,v}
&=(\Phi_{\la-\wt(v)-\wt(w)}^u\ten 1\ten 1)\circ\Phi_\la^{w,v}
\cr
&=P^{34}\sum_i(\Phi_{\la-\wt(v_i)-\wt(w_i)}^u\ten 1\ten 1)\circ
\Phi_\la^{v_i,w_i}
\cr
&=P^{34}\sum_i(\Phi_{\la-\wt(w_i)}^{u,v_i}\ten 1)\circ\Phi_\la^{w_i}
\cr
&=P^{34}P^{23}\sum_i\sum_j(\Phi_{\la-\wt(w_i)}^{(v_i)_j,u_j}\ten 1)
\circ\Phi_\la^{w_i}
\cr
&=P^{34}P^{23}\sum_i\sum_j(\Phi_{\la-\wt(w_i)-\wt(u_j)}^{(v_i)_j}\ten 1\ten 1)
\circ\Phi_\la^{u_j,w_i}
\cr
&=P^{34}P^{23}P^{34}\sum_i\sum_j\sum_k
(\Phi_{\la-\wt((w_i)_k)-\wt((u_j)_k)}^{(v_i)_j}\ten 1\ten 1)
\circ\Phi_\la^{(w_i)_k,(u_j)_k}
\cr
&=P^{34}P^{23}P^{34}\sum_i\sum_j\sum_k\Phi_\la^{(v_i)_j,(w_i)_k,(u_j)_k},
&\eq{2}
\cr}
$$
and accordingly
\goodbreak
$$
\eqalignno{
R_{WU}(\la)\,R_{VU}(\la-h^{(2)})\,R_{VW}(\la)\,(v\ten w\ten u)
&=\sum_i R_{WU}(\la)\,R_{VU}(\la-h^{(2)})\,(v_i\ten w_i\ten u)
\cr
&=\sum_i R_{WU}(\la)\,R_{VU}(\la-\wt(w_i))\,(v_i\ten w_i\ten u)
\cr
&=\sum_i\sum_j R_{WU}(\la)\,((v_i)_j\ten w_i\ten u_j)
\cr
&=\sum_i\sum_j\sum_k(v_i)_j\ten (w_i)_k\ten(u_j)_k.&\eq{3}
\cr}
$$

On the other hand we have
$$
\eqalignno{
\Phi_\la^{u,w,v}
&=(\Phi_{\la-\wt(v)}^{u,w}\ten 1)\circ\Phi_\la^v
\cr
&=P^{23}\sum_i(\Phi_{\la-\wt(v)}^{w_i,v_i}\ten 1)\circ\Phi_\la^v
\cr
&=P^{23}\sum_i(\Phi_{\la-\wt(v)-\wt(u_i)}^{w_i}\ten 1\ten 1)\circ
\Phi_\la^{u_i,v}
\cr
&=P^{23}P^{34}\sum_i\sum_j
(\Phi_{\la-\wt(v_j)-\wt((u_i)_j)}^{w_i}\ten 1\ten 1)
\circ\Phi_\la^{v_j,(u_i)_j}
\cr
&=P^{23}P^{34}\sum_i\sum_j
(\Phi_{\la-\wt((u_i)_j)}^{w_i,v_j}\ten 1)\circ\Phi_\la^{(u_i)_j}
\cr
&=P^{23}P^{34}P^{23}\sum_i\sum_j\sum_k
(\Phi_{\la-\wt((u_i)_j)}^{(v_j)_k,(w_i)_k}\ten 1)\circ\Phi_\la^{(u_i)_j}
\cr
&=P^{23}P^{34}P^{23}\sum_i\sum_j\sum_k\Phi_\la^{(v_j)_k,(w_i)_k,(u_i)_j},
&\eq{4}
\cr}
$$
and accordingly
$$
\eqalignno{
R_{VW}(\la-h^{(3)})\,&R_{VU}(\la)\,R_{WU}(\la-h^{(1)})\,(v\ten w\ten u)
\cr
&=R_{VW}(\la-h^{(3)})\,R_{VU}(\la)\,R_{WU}(\la-\wt(v))\,(v\ten w\ten u)
\cr
&=\sum_i R_{VW}(\la-h^{(3)})\,R_{VU}(\la)\,(v\ten w_i\ten u_i)
\cr
&=\sum_i\sum_j R_{VW}(\la-h^{(3)})\,(v_j\ten w_i\ten (u_i)_j)
\cr
&=\sum_i\sum_j R_{VW}(\la-\wt((u_i)_j))\,(v_j\ten w_i\ten (u_i)_j)
\cr
&=\sum_i\sum_j \sum_k (v_j)_k\ten (w_i)_k\ten(u_i)_j.&\eq{5}
\cr}
$$
It follows from \eqtag{2} and \eqtag{4} that
$$
\sum_i\sum_j\sum_k\Phi_\la^{(v_i)_j,(w_i)_k,(u_j)_k}=
\sum_i\sum_j\sum_k\Phi_\la^{(v_j)_k,(w_i)_k,(u_i)_j}.
$$
Hence the right-hand sides of \eqtag{3} and \eqtag{5} are equal.
Thus the left-hand sides of \eqtag{3} and \eqtag{5} are also equal.\hhalmos

\bPP
As pointed out in [E-S], \S2.2 the construction of intertwining operators,
fusion and exchange matrices admit natural quantum analogues.
Most definitions, results and proofs go on essentially unchanged compared
to the $q=1$ case. However, in the definition of the exchange matrix
the R-matrix $\FSR_{VW}$ associated to $U_q(\gog)$-modules $V$ and $W$,
and induced by the universal R-matrix $\FSR$,
is also needed. I will use the notation
$$
\FSR_{WV}^{21}:=(\FSR_{WV})^{21}=P_{WV} \FSR_{WV} P_{VW}.
\eqno\eq{35}
$$
This is different from the notation
$\FSR_{VW}^{21}:=(\FSR^{21})_{VW}$
in [E-S], \S2.2. The exchange matrix in the quantum case is now defined by
$$
R_{VW}(\la):=J_{VW}(\la)^{-1}\,\FSR_{WV}^{21}\,J_{WV}^{21}(\la).
\eqno\eq{36}
$$

The dynamical two-cocycle condition \eqtag{34} will remain valid in the
quantum case.
I will now discuss a second proof 
of the QDYBE \eqtag{39}, which is briefly sketched in
the remark in [E-S] after Proposition 2.4, and which also holds in the
quantum case. In the following, when being in the $q=1$ case, just put
$\FSR_{VW}$ equal to 1 (for any $V,W$).

I derive first the following
two important formulas (not given in [E-S]) for the exchange matrix:
$$
\eqalignno{
J_{VW}(\la)^{-1}\,R_{U,V\ten W}(\la)\,J_{VW}(\la-h^{(U)})
&=
R_{UV}(\la-h^{(W)})\,R_{UW}(\la),
&\eq{6}
\cr
J_{UV}(\la-h^{(W)})^{-1}\,R_{U\ten V,W}(\la)\,J_{UV}(\la)
&=
R_{VW}(\la)\,R_{UW}(\la-h^{(V)}),
&\eq{7}
}
$$
where both sides in \eqtag{6} and \eqtag{7} are acting on $U\ten V\ten W$.
Here we have adapted the notation introduced in [E-S] just before
Proposition 2.3 as follows. If $U=A_i$ then $F(\la-h^{(U)})$ will mean
$F(\la-h^{(i)})$.

One of the formulas \eqtag{6} and \eqtag{7} can be obtained by specialization
of formula (2.42) in [E-V]. Note that \eqtag{6} and \eqtag{7}
are also dynamical analogues
of the formulas
$$
\FSR_{U\ten V,W}=\FSR_{UW}\,\FSR_{VW},\quad
\FSR_{U,V\ten W}=\FSR_{UW}\,\FSR_{UV},
\eqno\eq{37}
$$
obtained from the following formulas for the universal R-matrix:
$$
(\De\ten \id)(\FSR)=\FSR_{13}\,\FSR_{23},\quad
(\id\ten\De)(\FSR)=\FSR_{13}\,\FSR_{12},
$$
which belong to the defining properties of a quasitriangular
Hopf algebra. Another defining property of a  quasitriangular
Hopf algebra is that
$$
P\,(\De(u))=\FSR\,\De(u)\,\FSR^{-1},
$$
which implies for the universal fusion matrix $J(\la)$ (see
[E-S], \S8) that
$$
\eqalignno{
&P_{12}\,(\De\ten 1)(J(\la))=\FSR_{12}\,(\De\ten\ 1)(J(\la))\,\FSR_{12}^{-1},
\cr
&P_{23}\,(1\ten\De)(J(\la))=\FSR_{23}\,(1\ten\De)(J(\la))\,\FSR_{23}^{-1},
}
$$
and hence
$$
\eqalign{
&P_{WV}\,J_{W\ten V,U}(\la)\,P_{VW}=
\FSR_{VW}\,J_{V\ten W, U}(\la)\,\FSR_{VW}^{-1}.
\cr
&P_{UW}\,J_{V,U\ten W}(\la)\,P_{WU}=
\FSR_{WU}\,J_{V,W\ten U}(\la)\,\FSR_{WU}^{-1}.
}
\eqno\eq{38}
$$
In the proof of \eqtag{6} and \eqtag{7} I will need
\eqtag{37} and \eqtag{38}.

\mLP
{\bf Proof of \eqtag{6}}
$$
\eqalignno{
&J_{VW}(\la)^{-1}\,R_{U,V\ten W}(\la)\,J_{VW}(\la-h^{(U)})
\cr
&=
J_{VW}(\la)^{-1}\,J_{U,V\ten W}(\la)^{-1}\,
J_{V\ten W,U}^{21}(\la)\,\FSR_{V\ten W,U}^{21}\,
J_{VW}(\la-h^{(U)})
\cr
&=
J_{UV}(\la-h^{(W)})^{-1}\,J_{U\ten V,W}(\la)^{-1}\,P_{VU}\,P_{WU}\,
\cr
&\qquad\qquad\qquad\qquad
\FSR_{V\ten W,U}\,J_{V\ten W,U}(\la)\,P_{UW}\,P_{UV}\,J_{VW}(\la-h^{(U)})
\cr
&=
J_{UV}(\la-h^{(W)})^{-1}\,J_{U\ten V,W}(\la)^{-1}\,P_{VU}\,P_{WU}
\cr
&\qquad\qquad\qquad\qquad
\FSR_{VU}\,\FSR_{WU}\,
J_{V\ten W,U}(\la)\,J_{VW}(\la-h^{(U)})\,P_{UW}\,P_{UV}
\cr
&=
J_{UV}(\la-h^{(W)})^{-1}\,J_{U\ten V,W}(\la)^{-1}\,P_{VU}\,
\FSR_{VU}\,P_{WU}\,\FSR_{WU}\,
J_{V,W\ten U}(\la)\,J_{WU}(\la)\,P_{UW}\,P_{UV}
\cr
&=
J_{UV}(\la-h^{(W)})^{-1}\,P_{VU}\,\FSR_{VU}\,J_{V\ten U,W}(\la)^{-1}\,
J_{V,U\ten W}(\la)\,P_{WU}\,\FSR_{WU}\,J_{WU}(\la)\,P_{UW}\,P_{UV}
\cr
&=
J_{UV}(\la-h^{(W)})^{-1}\,P_{VU}\,\FSR_{VU}\,J_{VU}(\la-h^{(W)})\,
J_{UW}(\la)^{-1}\,P_{WU}\,\FSR_{WU}\,J_{WU}(\la)\,P_{UW}\,P_{UV}
\cr
&=
J_{UV}(\la-h^{(W)})^{-1}\,\FSR_{VU}^{21}\,J_{VU}^{21}(\la-h^{(W)})\,P_{VU}\,
J_{UW}(\la)^{-1}\,\FSR_{WU}^{21}\,J_{WU}^{21}(\la)\,P_{UV}
\cr
&=
R_{UV}(\la-h^{(W)})\,R_{UW}(\la).&\halmos}
$$
\mLP
{\bf Proof of \eqtag{7}}
$$
\eqalignno{
&J_{UV}(\la-h^{(W)})^{-1}\,R_{U\ten V,W}(\la)\,J_{UV}(\la)
\cr
&=
J_{UV}(\la-h^{(W)})^{-1}\,J_{U\ten V,W}(\la)^{-1}\,
\FSR_{W,U\ten V}^{21}\,J_{W,U\ten V}^{21}(\la)\,J_{UV}(\la)
\cr
&=
J_{VW}(\la)^{-1}\,J_{U,V\ten W}(\la)^{-1}\,P_{WV}\,P_{WU}\,
\FSR_{W,U\ten V}\,J_{W,U\ten V}(\la)\,P_{UW}\,P_{VW}\,J_{UV}(\la)
\cr
&=
J_{VW}(\la)^{-1}\,J_{U,V\ten W}(\la)^{-1}\,P_{WV}\,P_{WU}\,
\FSR_{WV}\,\FSR_{WU}\,J_{W,U\ten V}(\la)\,J_{UV}(\la)\,P_{UW}\,P_{VW}
\cr
&=
J_{VW}(\la)^{-1}\,J_{U,V\ten W}(\la)^{-1}\,P_{WV}\,\FSR_{WV}\,P_{WU}\,
\FSR_{WU}\,J_{W\ten U,V}(\la)\,J_{WU}(\la-h^{(V)})\,P_{WU}\,P_{WV}
\cr
&=
J_{VW}(\la)^{-1}\,P_{WV}\,\FSR_{WV}\,J_{U,W\ten V}(\la)^{-1}\,
J_{U\ten W,V}(\la)\,P_{WU}\,\FSR_{WU}\,J_{WU}(\la-h^{(V)})\,P_{UW}\,P_{VW}
\cr
&=
J_{VW}(\la)^{-1}\,P_{WV}\,\FSR_{WV}\,J_{WV}(\la)\,
J_{UW}(\la-h^{(V)})^{-1}\,P_{WU}\,\FSR_{WU}\,
J_{WU}(\la-h^{(V)})\,P_{UW}\,P_{VW}
\cr
&=
J_{VW}(\la)^{-1}\,\FSR_{WV}^{21}\,J_{WV}^{21}(\la)\,P_{WV}\,
J_{UW}(\la-h^{(V)})^{-1}\,\FSR_{WU}^{21}\,J_{WU}^{21}(\la-h^{(V)})\,P_{VW}
\cr
&=
R_{VW}(\la)\,R_{UW}(\la-h^{(V)}).
&\halmos}
$$
In both proofs we have used the 2-cocycle condition \eqtag{34}
for the fusion matrix three times.
\bLP
{\bf Second proof of QDYBE \eqtag{39}}\quad (using \eqtag{6} and \eqtag{7};
acting on $V\ten W\ten U$)
$$
\eqalignno{
&R_{VW}(\la-h^{(U)})\,R_{VU}(\la)\,R_{WU}(\la-h^{(V)})
\cr
&=J_{WU}(\la)^{-1}\,R_{V,W\ten U}(\la)\,J_{WU}(\la-h^{(V)})\,
J_{WU}(\la-h^{(V)})^{-1}\,\FSR_{UW}^{21}\,J_{UW}^{21}(\la-h^{(V)})
\cr
&=J_{WU}(\la)^{-1}\,R_{V,W\ten U}(\la)\,P_{UW}\,\FSR_{UW}\,
J_{UW}(\la-h^{(V)})\,P_{WU}
\cr
&=J_{WU}(\la)^{-1}\,P_{UW}\,\FSR_{UW}\,
R_{V,U\ten W}(\la)\,J_{UW}(\la-h^{(V)})\,P_{WU}
\cr
&=J_{WU}(\la)^{-1}\,P_{UW}\,\FSR_{UW}\,J_{UW}(\la)\,R_{VU}(\la-h^{(W)})\,
R_{VW}(\la)\,P_{WU}
\cr
&=R_{WU}(\la)\,P_{UW}\,R_{VU}(\la-h^{(W)})\,R_{VW}(\la)\,P_{WU}
\cr
&=R_{WU}(\la)\,R_{VU}(\la-h^{(W)})\,R_{VW}(\la).
&\halmos}
$$

\Sec{difference} {Difference operators}
Next I give a proof for the $q=1$ case of the formula
$$
\FSD_{V\ten W}^U=\FSD_V^U\,\FSD_W^U=\FSD_W^U\,\FSD_V^U,
\eqno\eq{8}
$$
stated at the end of \S9.1 in [E-S] for the quantum case.
Let $\gog$ be a simple Lie algebra. For any two finite-dimensional
$\gog$-modules $U$ and $V$ let $R_{VU}(\la)$ be the exchange matrix.
Let $\RR_{VU}(\la):=R_{VU}(-\la-\rho)$ denote the shifted exchange matrix.
Let $\FSF_U$ be the space of $U[0]$-valued meromorphic functions on
$\goh^*$. For $\nu\in\goh^*$ let $T_\nu\in\End(\FSF_U)$ be the shift
operator $(T_\nu f)(\la):=f(\la+\nu)$.
Define the difference operator $\FSD_V^{\la,U}$ acting on $\FSF_U$ by
$$
\eqalignno{
\FSD_V^{\la,U}:=&\sum_{\nu\in\goh^*}\Tr|_{V[\nu]}\,(\RR_{VU}(\la))\,T_\nu
\cr
=&\sum_{\nu\in\goh^*}\Tr|_{V[\nu]}\,(\RR_{V[\nu],U[0];V[\nu],U[0]}(\la))\,
T_\nu,&\eq{31}
}
$$
where $\RR_{V[\la],U[\mu];V[\nu],U[\si]}$ denotes the block of the matrix
$\RR_{VU}$ corresponding to the weight spaces
$V[\la],U[\mu];V[\nu],U[\si]$ (which block will be zero unless
$\la+\mu=\nu+\si$).

\bLP
{\bf Proof of \eqtag{8}}
\sLP
We can rewrite \eqtag{7} as
$$
\RR_{W\ten V,U}(\la)=\JJ_{WV}(\la+h^{(U)})\,\RR_{VU}(\la)\,
\RR_{WU}(\la+h^{(V)})\,\JJ_{WV}(\la)^{-1},
$$
where $\JJ_{WV}(\la):=J(-\la-\rho)$ denotes the shifted fusion matrix.
Hence
$$
\eqalignno{
&\RR_{W[\nu]\ten V[\mu],U[0];W[\nu]\ten V[\mu],U[0]}(\la)=
\sum_{\mu',\nu',\mu'',\nu'',\si}\JJ_{W[\nu],V[\mu];W[\nu'],V[\mu']}(\la)\,
\cr
&\quad\circ \RR_{V[\mu'],U[0];V[\mu''],U[\si]}(\la)\,
\RR_{W[\nu'],U[\si];W[\nu''],U[0]}(\la+\mu'')\,
\JJ_{W[\nu],V[\mu];W[\nu''],V[\mu'']}(\la)^{-1}.
}
$$
Hence
$$
\eqalignno{
&\Tr|_{W[\nu]\ten V[\mu]}\,
(\RR_{W[\nu]\ten V[\mu],U[0];W[\nu]\ten V[\mu],U[0]}(\la))
\cr
&\quad=
\sum_\si\Tr|_{W[\nu]\ten V[\mu]}\,
(\RR_{V[\mu],U[0];V[\mu],U[\si]}(\la)\,
\RR_{W[\nu],U[\si];W[\nu],U[0]}(\la+\mu))
\cr
&\quad=
\Tr|_{W[\nu]\ten V[\mu]}\,
(\RR_{V[\mu],U[0];V[\mu],U[0]}(\la)\,
\RR_{W[\nu],U[0];W[\nu],U[0]}(\la+\mu)).
}
$$
Then
$$
\eqalignno{
\FSD_V^{\la,U}\,\FSD_W^{\la,U}
&=\sum_\mu\Tr|_{V[\mu]}\,(\RR_{V[\mu],U[0];V[\mu],U[0]}(\la))\,
T_\mu\,\sum_\nu\Tr|_{W[\nu]}\,(\RR_{W[\nu],U[0];W[\nu],U[0]}(\la))\,T_\nu
\cr
&=\sum_{\mu,\nu}\Tr|_{V[\mu]}\,(\RR_{V[\mu],U[0];V[\mu],U[0]}(\la))\,
\Tr|_{W[\nu]}\,(\RR_{W[\nu],U[0];W[\nu],U[0]}(\la+\mu))\,T_{\mu+\nu}
\cr
&=\sum_{\mu,\nu}\Tr|_{W[\nu]\ten V[\mu]}\,
(\RR_{V[\mu],U[0];V[\mu],U[0]}(\la)\,
\RR_{W[\nu],U[0];W[\nu],U[0]}(\la+\mu))\,T_{\mu+\nu}
\cr
&=\sum_{\mu,\nu}
\Tr|_{W[\nu]\ten V[\mu]}\,
(\RR_{W[\nu]\ten V[\mu],U[0];W[\nu]\ten V[\mu],U[0]}(\la))\,T_{\mu+\nu}
\cr
&=\sum_\si\Tr|_{(W\ten V)[\si]}\,
(\RR_{(W\ten V)[\si],U[0];(W\ten V)[\si],U[0]}(\la))\,T_\si
=\FSD_{W\ten V}^{\la,U}.
}
$$
But also
$$
\eqalignno{
&\;\FSD_{W\ten V}^{\la,U}
=\sum_{\mu,\nu}
\Tr|_{W[\nu]\ten V[\mu]}\,
(\RR_{W[\nu]\ten V[\mu],U[0];W[\nu]\ten V[\mu],U[0]}(\la))\,T_{\mu+\nu}
\cr
&=\sum_{\mu,\nu}
\Tr|_{W[\nu]\ten V[\mu]}\,(P_{V[\mu],W[\nu];V[\mu],W[\nu]}\,
\RR_{V[\mu]\ten W[\nu],U[0];V[\mu]\ten W[\nu],U[0]}(\la)\,
P_{W[\nu],V[\mu];W[\nu],V[\mu]})
\cr
&\qquad\circ T_{\mu+\nu}=\sum_{\mu,\nu}
\Tr|_{V[\mu]\ten W[\nu]}\,
(\RR_{V[\mu]\ten W[\nu],U[0];V[\mu]\ten W[\nu],U[0]}(\la))\,T_{\mu+\nu}
=\FSD_{V\ten W}^{\la,U}\,.
}
$$
Hence
$$
\FSD_V^{\la,U}\,\FSD_W^{\la,U}=\FSD_{W\ten V}^{\la,U}=\FSD_{V\ten W}^{\la,U}
=\FSD_W^{\la,U}\,\FSD_V^{\la,U}.
\eqno\halmos
$$

Note that it was possible to apply \eqtag{7} in the above proof above
because we had
assumed thet $\FSD_V^{\la,U}$ acts on $U[0]$-valued functions, and because the
definition of $\FSD_V^{\la,U}$ involved shift operators $T_\nu$.

\Sec{weightedtrace} {Weighted trace functions}
In \S9.2 of [E-S] weighted-trace functions are introduced and difference
equations are given for them. [E-S] refers for the proofs to [E-V].
Theorem 9.2 of [E-S] survives for $q=1$, see [E-V], Theorem 10.4.
I will give a proof of that theorem parallel to the proof of the $q$-case,
see Theorem 1.2 and \S3 in [E-V].

First consider the proof of Lemma 2.14 in [E-V].
Let $W$ be a finite-dimensional $\gog$-module. By the properties of
the intertwining operator we can uniquely define a bilinear form
$B_{\la,W}\colon W\times W^*\to\CC$ by the formula
$$
(1\ten\lan,\ran)\circ(\Phi_{\la-\wt(w^*)}^w\ten 1)\circ\Phi_\la^{w^*}=
B_{\la,W}(w,w^*)\,\id_{M_\la}.
\eqno\eq{15}
$$
Note that $B_{\la,W}(w,w^*)=0$ if
$\wt(w)+\wt(w^*)\ne0$.
Since
$$
\Phi_\la^{J_{WW^*}(\la)(w\ten w^*)}=
(\Phi_{\la-\wt(w^*)}^w\ten 1)\circ\Phi_\la^{w^*},
$$
we have
$$
B_{\la,W}(w,w^*)=\lan,\ran\,\bigl( J_{WW^*}(\la)(w\ten w^*)\bigr).
\eqno\eq{16}
$$
Define a generalized element $Q(\la)$ in $U(\gog)$ in terms of the universal
fusion matrix by
$$
Q(\la):=\bigl(m\circ P\circ(1\ten S^{-1})\bigr)\,J(\la).
\eqno\eq{14}
$$
This induces an endomorphism $Q_W(\la)$ of $W$ given by
$$
Q_W(\la)=C_W\bigl( (J_{WW^*}(\la)^{t_2})^{21}\bigr),
\eqno\eq{13}
$$
where $C_W$ denotes contraction of an endomorphism of $W\ten W$ to
an endomorphism of $W$.
Now we have
$$
B_{\la,W}(w,w^*)=\lan Q_W(\la)\,w,w^*\ran.
\eqno\eq{12}
$$
Indeed, if $T\in\End(W\ten W)$ then
$\lan C(T)w,w^*\ran=\lan,\ran\bigl((T^{21})^{t_2}(w\ten w^*)\bigr)$. Hence
$$
\lan Q_W(\la)\,w,w^*\ran=
\lan,\ran\bigl(J_{WW^*}(\la)(w\ten w^*)\bigr)=B_{\la,W}(w,w^*).
$$
It follows from \eqtag{12} that $Q_W(\la)$ is a weight preserving
endomorphism of $W$.

Next we have
$$
B_{\la,U\ten W}\circ\bigl(J_{UW}(\la-h^{(U^*)}-h^{(W^*)})\ten
J_{U^*W^*}(\la)\bigr)=B_{\la,U}\circ B_{\la-h^{(U^*)},W}.
\eqno\eq{17}
$$
Indeed,
$$
\eqalignno{
&B_{\la,U}(u,u^*)\,B_{\la-\wt(u^*),W}\,\id_{M_\la}
\cr
&=
\bigl(\lan,\ran\ten\lan,\ran\bigr)\circ
\Phi_{\la-\wt(u^*)-\wt(w^*)-\wt(w)}^u\circ\Phi_{\la-\wt(u^*)-\wt(w^*)}^w\circ
\Phi_{\la-\wt(u^*)}^{w^*}\circ\Phi_\la^{u^*}
\cr
&=\bigl(\lan,\ran\ten\lan,\ran\bigr)\circ
\Phi_{\la-\wt(u^*)-\wt(w^*)}^{J_{UW}(\la-\wt(u^*)-\wt(w^*))(u\ten w)}\circ
\Phi_\la^{J_{W^*U^*}(\la)(w^*\ten u^*)}
\cr
&=B_{\la,U\ten W}\bigl(J_{UW}(\la-\wt(u^*)-\wt(w^*))(u\ten w),
J_{W^*U^*}^{21}(\la)(u^*\ten w^*)\bigr)\,\id_{M_\la}.
}
$$
Combination of \eqtag{17} with \eqtag{12} yields that
$$
Q_{U\ten W}(\la)=\bigl(J_{W^*U^*}^{t_1t_2,21}(\la)\bigr)^{-1}\,
\bigl(Q_U(\la)\ten Q_W(\la+h^{(U)})\bigr)\,
J_{UW}(\la+h^{(U)}+h^{(W)})^{-1}.
\eqno\eq{18}
$$
It follows from \eqtag{12} and \eqtag{15} that
$Q_{U\ten W}(\la)=Q_{W\ten U}^{21}(\la)$. Hence
we can rewrite \eqtag{18} as
$$
Q_{U\ten W}(\la)=\bigl(J_{U^*W^*}^{t_1t_2}(\la)\bigr)^{-1}\,
\bigl(Q_U(\la+h^{(W)})\ten Q_W(\la)\bigr)\,
J_{WU}^{21}(\la+h^{(U)}+h^{(W)})^{-1}.
$$
Now eliminate $Q_{U\ten W}(\la)$ from these two formulas and substitute
$$
R_{UW}(\la)=J_{UW}(\la)^{-1}\,J_{WU}^{21}(\la)
\eqno\eq{19}
$$
(the defining formula for the exchange matrix in \S2.1 of [E-S]).
Then we obtain
$$
\eqalignno{
R_{U^*W^*}^{t_1t_2}(\la)=&
\bigl(Q_U(\la)\ten Q_W(\la+h^{(U)})\bigr)
\cr
&\qquad\circ
R_{UW}(\la+h^{(U)}+h^{(W)})\,
\bigl(Q_U(\la+h^{(W)})\ten Q_W(\la)\bigr)^{-1}.&\eq{20}
}
$$
This is essentially the formula at the end of \S3.3 in [E-V].

\bPP
Next I discuss Proposition 3.1 in [E-V].
Fix finite-dimensional $\gog$-modules $V$ and $W$. Let $B$ be a basis of $V$
consisting of weight vectors. For $v\in B$ let $V^*$ be the corresponding
dual basis vector of $V^*$. Define the operator
$$
\Phi_\mu^V\colon\quad y\mapsto \sum_{v\in B}\Phi_\mu^v(y)\ten v^*\colon\quad
M_\mu\to\bigoplus_\lambda\bigl(M_{\mu-\la}\ten V\ten V^*[-\la]\bigr),
\eqno\eq{23}
$$
which is clearly independent of the choice of $B$.
Define the isomorphism
$$
\eta_W(\mu)\colon\quad \bigoplus_\nu\bigl(W[\nu]\otimes M_{\mu+\nu}\bigr)\to
M_\mu\otimes W,
$$
where
$$
\eta_W(\mu)(w\ten z):=\Phi_{\mu+\nu}^w(z)\quad
\hbox{if\quad $w\in W[\nu]$, $z\in M_{\mu+\nu}$.}
\eqno\eq{21}
$$
Proposition 3.1 together with formula (3.2) in [E-V] can now be formulated
as follows:
$$
\eqalignno{
&P_{V\ten V^*,W}\circ(\Phi_\mu^V\ten\id_W)\circ\eta_W(\mu)
\bigr|_{W[\nu]\ten M_{\mu+\nu}}
\cr
&\qquad=
(\eta_W(\mu)
\ten\id_V\ten\id_{V^*})\circ R_{WV}^{t_2}(\mu+\nu)\circ
(\id_W\ten\Phi_{\mu+\nu}^V)\bigr|_{W[\nu]\ten M_{\mu+\nu}}.
&\eq{22}
}
$$

\mLP
{\bf Proof of \eqtag{22}}\quad
Write
$R_{WV}(\mu+\nu)=\sum_i p_i\ten q_i^t$.
Let $y\in M_{\mu+\nu}$ and $w\in W[\nu]$. Then
$$
\eqalignno{
&\bigl(P_{V\ten V^*,W}\circ(\Phi_\mu^V\ten\id_W)\circ\eta_W(\mu)\bigr)(w\ten y)
=\bigl(P_{V\ten V^*,W}\circ(\Phi_\mu^V\ten\id_W)\circ\Phi_{\mu+\nu}^w\bigr)(y)
\cr
&\quad
=P_{VW}\sum_{v\in B}(\Phi_\mu^v\ten\id)\bigl(\Phi_{\mu+\nu}^w(y)\bigr)\ten v^*
=P_{VW}\sum_{v\in B}\Phi_{\mu+\nu}^{J_{VW}(\mu+\nu)(v\ten w)}(y)\ten v^*
\cr
&\quad
=\sum_{v\in B}\Phi_{\mu+\nu}^{J_{VW}^{21}(\mu+\nu)(w\ten v)}(y)\ten v^*
=\sum_{v\in B}\Phi_{\mu+\nu}^{J_{WV}(\mu+\nu)(R_{WV}(\mu+\nu)(w\ten v))}
(y)\ten v^*
\cr
&\qquad
=\sum_{v\in B}\sum_i(\Phi_{\mu+\nu-\wt(v)}^{p_iw}\ten\id_V)
\bigl(\Phi_{\mu+\nu}^{q_i^tv}(y)\bigr)\ten v^*
\cr
&\qquad
=\sum_{v\in B}\sum_i(\Phi_{\mu+\nu-\wt(v)}^{p_iw}\ten\id_V)
\bigl(\Phi_{\mu+\nu}^v(y)\bigr)\ten q_i v^*
\cr
&\qquad=\sum_{v\in B}\sum_i(\eta_W(\mu)\ten\id_V\ten\id_{V^*})
(p_iw\ten\Phi_{\mu+\nu}^v(y)\ten q_i v^*)
\cr
&\qquad=(\eta_W(\mu)
\ten\id_V\ten\id_{V^*})\circ R_{WV}^{t_2}(\mu+\nu)\circ
(w\ten\Phi_{\mu+\nu}^V(y)).&\halmos
}
$$

\bPP
For $\la\in\goh^*$ and $U$ a $\gog$-module let
$e^\la\colon u\mapsto e^{\lan\la,\wt(u)\ran}\,u\colon U\to U$.
Let $V$ be a finite dimensional $\gog$-module and let $v\in V[0]$.
Let $\{y_i\}$ be a basis of $M_\mu$ consisting of weight vectors.
Since $\Phi_\mu^v\colon M_\mu\to M_\mu\ten V$ is
weight preserving, we have
$\Phi_\mu^v(y_i)\in y_i\ten V[0]+\sum_{j\ne i} y_j\ten V$.
Hence, if $B[0]$ is a basis of $V[0]$ and
if $v^*\in V^*[0]$ is the dual basis vector corresponding to $v\in B[0]$,
we have
$$
\Phi_\mu^v(e^\la y_i)\ten v^*\in y_i\ten V[0]\ten V^*[0]+
\sum_{j\ne i} y_j\ten V\ten V^*[0].
$$
Let
$$
\Phi_\mu^{V[0]}:=\sum_{v\in B[0]}\Phi_\mu^v\ten v^*.
\eqno\eq{24}
$$
Then
$$
\Psi_V(\la,\mu):=\Tr\bigl|_{M_\mu}(\Phi_\mu^{V[0]}\circ
e^\la)\in V[0]\ten V^*[0].
\eqno\eq{25}
$$
For $W$ a finite dimensional $\gog$-module let
$$
\chi_W(e^\la):=\Tr\bigl|_W e^\la=\sum_\nu \dim(W[\nu])\,e^{\lan\la,\nu\ran}.
\eqno\eq{26}
$$

A difference equation for $\Psi_V(\la,\mu)$ in the variable $\mu$ can be
derived from \eqtag{22}. First multiply both sides of
\eqtag{22} with $e^{\lan \la,\mu+\nu\ran}$, observe that
$\eta_W(\mu)\circ(\id_W\ten e^\la)=(e^\la\ten e^\la)\circ \eta_W(\mu)$,
sum both sides of \eqtag{22} with respect to $\nu$, and
next multiply both sides
of \eqtag{22} on the left with $(\eta_W(\mu)\ten\id_V\ten\id_{V^*})^{-1}$.
Then we obtain the following identity of linear endomorphisms
$\dirsum_\nu(W[\nu]\ten M_{\mu+\nu})\to
\dirsum_\nu(W[\nu]\ten M_{\mu+\nu}\ten V\ten V^*)$.
$$
\eqalignno{
&(\eta_W(\mu)\ten\id_V\ten\id_{V^*})^{-1}\circ
P_{V\ten V^*,W}\circ\bigl((\Phi_\mu^V\circ e^\la)\ten (e^\la\circ\id_W\bigr)
\circ\eta_W(\mu)
\bigr|_{\dirsum_\nu(W[\nu]\ten M_{\mu+\nu})}
\cr
&\qquad\qquad\qquad\quad=
\sum_\nu R_{WV}^{t_2}(\mu+\nu)\circ
\bigl((\id_W\ten(\Phi_{\mu+\nu}^V\circ e^\la)\bigr)
\bigr|_{W[\nu]\ten M_{\mu+\nu}}.
&\eq{27}
}
$$
Now take the trace with respect to $\dirsum_\nu(W[\nu]\ten M_{\mu+\nu})$
on both sides of \eqtag{27} and use \eqtag{24}. Then
$$
\bigl(\Tr\bigl|_W e^\la\bigr)\,
\bigl(\Tr\bigl|_{M_\mu}\Phi_\mu^{V[0]}\circ e^\la\bigr)=
\sum_\nu\Tr\bigl|_{W[\nu]} R_{WV}^{t_2}(\mu+\nu)\circ
\bigl(\id_W\ten\Tr\bigl|_{M_{\mu+\nu}}
\Phi_{\mu+\nu}^{V[0]}\circ e^\la\bigr).
$$
Now substitute \eqtag{25} and \eqtag{26}, and take inside the sum on the
right-hand side the transpose with respect to $W$. Then
$$
\chi_W(e^\la)\,\Psi_V(\la,\mu)=\sum_\nu\Tr\bigl|_{W^*[-\nu]}
R_{WV}^{t_1t_2}(\mu+\nu)\circ
\bigl(\id_{W^*}\ten\Psi_V(\la,\mu+\nu)\bigr).
\eqno\eq{28}
$$
On the right-hand side of formula \eqtag{28} substitute \eqtag{20}.
Next also substitute
$\RR_{VU}(\la):=R_{VU}(-\la-\rho)$ and $\QQ_{V}(\la):=Q_{V}(-\la-\rho)$.
Then
$$
\eqalignno{
&\chi_W(e^\la)\,\Psi_V(\la,\mu)
\cr
&=\sum_\nu\Tr\bigl|_{W^*[-\nu]}
\bigl(Q_{W^*}(\mu+\nu)\ten Q_{V^*}(\mu+\nu+h^{(W^*)})\bigr)
\circ
R_{W^*V^*}(\mu+\nu+h^{(V^*)}+h^{(W^*)})
\cr
&\qquad\qquad\circ
\bigl(Q_{W^*}(\mu+\nu+h^{(V^*)})\ten Q_{V^*}(\mu+\nu)\bigr)^{-1}
\circ\bigl(\id_{W^*}\ten\Psi_V(\la,\mu+\nu)\bigr)
\cr
&=
\sum_\nu\Tr\bigl|_{W^*[\nu]}
\bigl(\QQ_{W^*}(\mu+\nu)\ten \QQ_{V^*}(\mu+\nu-h^{(W^*)})\bigr)
\circ
\RR_{W^*V^*}(\mu+\nu-h^{(V^*)}-h^{(W^*)})
\cr
&\qquad\qquad\circ
\bigl(\QQ_{W^*}(\mu+\nu-h^{(V^*)})\ten \QQ_{V^*}(\mu+\nu)\bigr)^{-1}
\circ\bigl(\id_{W^*}\ten\Psi_V(\la,-\mu-\nu-\rho)\bigr)
\cr
&=
\sum_\nu\Tr\bigl|_{W^*[\nu]}
\bigl(\QQ_{W^*}(\mu+\nu)\ten \QQ_{V^*}(\mu)\bigr)
\circ
\RR_{W^*V^*}(\mu)\circ
\bigl(\QQ_{W^*}(\mu+\nu)\ten \QQ_{V^*}(\mu+\nu)\bigr)^{-1}
\cr
&\qquad\qquad\circ\bigl(\id_{W^*}\ten\Psi_V(\la,-\mu-\nu-\rho)\bigr)
\cr
&=
\QQ_{V^*}(\mu)\circ\sum_\nu\Tr\bigl|_{W^*[\nu]}
\RR_{W^*V^*}(\mu)
\circ\Bigl(\id_{W^*}\ten
\bigl(\QQ_{V^*}(\mu+\nu)^{-1}\circ\Psi_V(\la,-\mu-\nu-\rho)\bigr)\Bigr).
\cr
&&\eq{32}
}
$$
Let the Weyl denominator be given by
$$
\de(\la):=e^{\lan\la,\rho\ran}\prod_{\al>0}(1-e^{-\lan\la,\al\ran}).
\eqno\eq{30}
$$
Define the weighted-trace function by
$$
F_V(\la,\mu):=\QQ_{V^*}^{-1}(\mu)\,\Psi_V(\la,-\mu-\rho)\,\de(\la).
\eqno\eq{29}
$$
Now replace $W^*$ by $W$ in \eqtag{32} and substitute \eqtag{31} and \eqtag{29}
in \eqtag{32}. We finally obtain the formula which is the $q=1$ case of
Theorem 9.2 in [E-S]:
$$
\FSD_W^{\mu,V^*}\,F_V(\la,\mu)=\chi_W(e^{-\la})\,F_V(\la,\mu).
\eqno\eq{33}
$$

\Ref
\refitem B-B-B
O. Babelon, D. Bernard \& E. Billey,
{\sl A quasi-Hopf algebra interpretation of quantum 3-j and 6-j symbols
and difference equations},
Phys. Lett. B 375 (1996), 89--97;
{\tt \hbox{q-alg/9511019}}.

\refitem E-S
P. Etingof \& O. Schiffmann,
{\sl Lectures on the dynamical Yang-Baxter equations},
preprint {\tt math.QA/9908064 v2}, 1999, 2000.

\refitem E-V1
P. Etingof \& A. Varchenko,
{\sl Geometry and classification of solutions of the classical dynamical
Yang-Baxter equation},
Comm. Math. Phys. 192 (1998), 77--120;
\item{}
{\tt q-alg/9703040}.

\refitem E-V2
P. Etingof \& A. Varchenko,
{\sl Solutions of the quantum dynamical Yang-Baxter equation and dynamical
quantum groups},
Comm. Math. Phys. 196 (1998), 591--640;
\item{} {\tt q-alg/9708015}.

\refitem E-V3
P. Etingof \& A. Varchenko,
{\sl Exchange dynamical quantum groups},
Comm. Math. Phys. 205 (1999), 19--52; {\tt math.QA/9801135}.

\refitem E-V4
P. Etingof \& A. Varchenko,
{\sl Traces of intertwiners for quantum groups and difference equations, I},
{\tt math.QA/9907181 v2}, 1999, 2000.

\refitem Fe1
G. Felder,
{\sl Conformal field theory and integrable systems associated to elliptic
curves},
in {\sl Proceedings of the International Congress of Mathematicians,
Z\"urich, 1994},
\item{}
Birkh\"auser 1994, pp. 1247--1255.

\refitem Fe2
G. Felder,
{\sl Elliptic quantum groups},
in
{\sl XIth International Congress of Mathematical Physics (Paris, 1994)},
Internat. Press, Cambridge, MA, 1995, pp. 211--218;
\item{} {\tt hep-th/9412207}.

\refitem G-N
J.-L. Gervais and A. Neveu,
{\sl Novel triangle relations and absence of tachyons in Liouville
string field theory},
Nucl. Phys. B 238 (1984), 125--141.

\nonfrenchspacing
\vskip 1truecm
\LP
{\obeylines\parindent 7truecm
\addressfont
Tom H. Koornwinder,
Korteweg-de Vries Institute, Universiteit van Amsterdam,
\LP
Plantage Muidergracht 24, 1018 TV Amsterdam,
The Netherlands;
\LP
email: {\ttaddressfont thk@science.uva.nl}}

\bye